\documentclass[12pt, reqno]{amsart}
\usepackage{amscd,times,hyperref}
\usepackage[enableskew,vcentermath]{youngtab}
\usepackage{amssymb}
\usepackage{amsmath}
\usepackage[all]{xypic}
\usepackage{graphics}

\newtheorem{thm}{Theorem}[section] 
\newtheorem{lem}[thm]{Lemma}

\newtheorem{cor}[thm]{Corollary}
\newtheorem{pro}[thm]{Proposition}
\newcommand{\bbar}{\begin{array}}
\newcommand{\eear}{\end{array}}
\newcommand{\bb}{\begin{equation}}
\newcommand{\eqbb}{\begin{equation}} 
\def\det{\mbox{\rm det}}
\def\ch{\mbox{\rm ch}}
 
\begin{document} 
 \bibliographystyle{plain}
\title[Jacobi-Trudi type formula]{Jacobi-Trudi type formula for character of  irreducible representations of  $\frak{gl}(m|1)$}
\author{Nguy\^en Luong Th\'ai B\`inh}
\address[NLT Binh]{Sai Gon University, Ho Chi Minh City, Vietnam}
 \email{nltbinh@sgu.edu.vn}

\author{Nguy\^en Thi Phuong Dung}
\address[NTP Dung]{Banking Academy, Hanoi, Vietnam} 
\email{dungnp@hvnh.edu.vn}

\author{Ph\`ung H\^o Hai}
\address[PH Hai]{Institute of Mathematics, Vietnam Academy of Science and Technology, Hanoi, Vietnam}
 \email{phung@math.ac.vn}

\dedicatory{Dedicated to Professor L\^e Tu\^an Hoa on the occasion of this sixtieth birthday}
 
\maketitle 
\begin{abstract}
We prove a determinantal type formula to compute the irreducible characters of the general Lie superalgebra $\mathfrak{gl}(m|1)$ in terms of the characters of the symmetric powers of the fundamental representation and their duals. This formula was conjectured by J. van der Jeugt and E. Moens for the Lie superalgebra $\frak{gl}(m|n)$ and generalizes the well-known Jacobi-Trudi formula.
\end{abstract}
\section{Introduction}
The classical Jacobi-Trudi formula computes Schur symmetric functions in terms of the elementary (resp. complete) symmetric functions. Since these symmetric functions can be realized as irreducible characters of general linear Lie algebras, we can interpret the Jacobi-Trudi formula as a formula for computing irreducible characters of general linear Lie algebras in terms of the characters of symmetric (resp. anti-symmetric) tensor representations. This formula complements the Weyl determinantal formula which computes irreducible characters in terms of root systems. Although the Jacobi-Trudi formula is well-defined only for partitions, that is, for integral dominant weight with non-negative components, it is well-known that an integral dominant weight can be led to a partition by adding some multiple of the partition $(1,1,\ldots,1)$, which corresponds to the  determinantal representation.

The aim of this work is to extend this famous formula to the case of the general linear Lie superalgebras. According to V.~Kac, irreducible (finite dimensional) representations of the general linear Lie superalgebra $\mathfrak{gl}(m|n)$ are determined by means of dominant weights. We shall restrict ourselves to representations with integral dominant weights, see Eq. \eqref{eq:lamda1}. V.~Kac established in the 70s an analog of Weyl formula to compute {  typical} irreducible characters. It took twenty years until V. Serganova provided a method to compute {\em atypical} irreducible
characters, which was subsequently simplified by Brundan \cite{Brundan} and Su-Zhang \cite{Zhang2}.

For some classes of integral dominant weights, Jacobi-Trudi formula has been established, for instance, when the weights correspond to partitions, i.e. the corresponding representation is constructed from the fundamental representation using multi-linear algebra. However, due to the more complicated nature of the representation categories of general linear Lie superalgebras, to extend Jacobi-Trudi formula to characters of mixed representations we will need to incorporate characters of both symmetric tensor powers and their duals. A conjectural determinantal formula was presented in detail in \cite{M6}. In fact, there was an unsuccessful attempt to prove it in \cite{M4} (see Concluding remarks for more details). 

In this work we prove the above mentioned determinantal formula for the case of $\mathfrak{gl}(m|1)$ (Theorem \ref{T:1}). In fact, we prove the formula for a class of integral dominant weights, which are called special weights (Equation \eqref{eq:15}).  As in the classical case, we notice that any integral dominant weight differs from a unique special weight by a multiple of the weight $(1,1,\ldots,1;-1)$, which corresponds to the super-determinantal representation (Proposition  \ref{L:3}).\\

\noindent{\bf In this work we shall agree with the following notations}
\begin{enumerate}
\item A partition is a finite decreasing sequence of non-negative integers. Given a partition $\lambda$, 
its length $l(\lambda)$ is the number of its positive components, 
its contents (or weight) $|\lambda|$ is the sum of its components. 
\item A dominant integral weight (of size $m$) is an ordered $m$-tuple of decreasing integers. Given a weight $\lambda$, its contents (or weight) is the sum of its components. A partition can be considered as a dominant weight in the obvious manner.\item The opposite to a partition $\lambda=(\lambda_1,\ldots,\lambda_p)$ with length $p\leq m$ is the weight
$(-\lambda_p,\ldots,-\lambda_1,0,\ldots, 0)$. It is denoted by $\overline\lambda$.

\item There is a natural operation of (component wise) addition  on the set of partitions (resp. dominant integral weights).
\item A composite partition is just a pair of partitions, it is called $m$-standard if the total length of the two partitions does not exceed $m$.
\item A integral dominant (super-) weight of size $(m|n)$ is an ordered $m+n$-tuple of integers such that the first $m$ components and the last $n$ components form non-increasing sequences.

\end{enumerate}
\section{Preliminaries}
This section presents some results on the general linear Lie superalgebras for the later use. We shall work over the complex field $\mathbb{C}$. \\

A vector super-space is a $\mathbb{Z/{\rm 2}Z}$-graded vector space $V = V_{\bar{0}}\oplus V_{\bar{1}}$. The vector spaces $ V_{\bar{0}}, V_{\bar{1}}$ are called even and odd homogeneous components of $V$, their elements are also called homogeneous elements.
A homogeneous element $x \in V_{\bar{0}}$ has degree   $\deg(x) = \bar 0$, while $x \in V_{\bar{1}}$ has degree  $\deg(x) = 1.$\\

Let $\text{End}(V)$ be a set of linear endomorphisms of $V$ then $\text{End}(V) = \text{End}_{\bar{0}}(V) \oplus \text{End}_{\bar{1}}(V)$, where
\begin{equation}\label{eq:}
\text{End}_{\bar{0}}(V) = \text{End}(V_{\bar{0}}) \oplus \text{End}(V_{\bar{1}})\;\; \mbox{and}\;\; \text{End}_{\bar{1}}(V) =\text{Hom}(V_{\bar{0}}, V_{\bar{1}}) \oplus \text{Hom}(V_{\bar{1}}, V_{\bar{0}}).
\end{equation}
We can equip $\text{End}(V)$ with the structure of a Lie superalgebra by defining the Lie bracket  $[-,-]$
\begin{equation}\label{eq:}
[x,y] = xy - (-1)^{\deg(x)\deg(y)}yx,
\end{equation}
on homogeneous elements and then extending linearly to all of $\text{End}(V)$.
Denote by $\frak{gl}(m|n)$ the vector super-space $\text{End}(V)$ equipped the above bracket, with $V = V_{\bar{0}}\oplus V_{\bar{1}}$, where $\dim V_{\bar{0}} = m, \dim V_{\bar{1}} = n$. In this paper we shall focus ourselves on $\mathfrak{gl}(m|1)$ .

\subsection{The Lie superalgebra $\frak{gl}(m|1)$}
In this paper we shall focus ourselves on $\frak{gl}(m|1)$, which is denoted by $\frak{g}$. We realize $\frak{g}$ as the set of $ (m + 1) \times (m + 1)$ matrices.
Hence 
\begin{equation}\label{eq:}
\frak{g}_{\bar{0}} = \left\{ \left (\begin{array}{cc}A&0\\
0&D \end{array}\right) | A \in M_{m,m}, D \in M_{1,1} \right\} \end{equation}
and
 \begin{equation} \frak{g}_{\bar{1}} = \left\{ \left (\begin{array}{cc}0&B\\
C&0 \end{array}\right) | B \in M_{m,1}, C \in M_{1,m} \right\},
\end{equation}
here $M_{r,t}$ denotes the set of $r\times t$ matrices.\\

The standard basis for $\frak{g}$ consists of matrices $E_{i,j}: i,j = 1, 2, \ldots, m + 1$ with 1 on the entry $(i, j)$ and $0$ elsewhere. Consider the subalgebra $\frak h \subset  \frak g$  spanned by the elements $ E_{j,j}: j = 1, 2, \ldots, m + 1$, $\frak{h}$ is a Cartan subalgebra of $\frak{g}$. The dual vector space  $\frak{h}^*$  is spanned by $\{\epsilon_i, \delta_1| i = 1, 2, \ldots, m  \}$, where $\epsilon_i(E_{j,j}) = \delta_{ij}$ and $\delta_1(E_{j,j}) = - \delta_{(m + 1)j}$ .
The roots of $\frak g$ can be expressed in terms of this basis. In this so-called distinguished choice for a triangular decomposition of $\frak g$, the simple roots are\\
$$\Pi = \{\epsilon_1-\epsilon_2, \cdots, \epsilon_{m-1}-\epsilon_m, \epsilon_m-\delta_1 \}.$$
In that case, the positive even roots are   
$$\Delta_{0}^+= \{\epsilon_i - \epsilon_j| 1 \leq i <  j \leq m\},$$
 and the positive odd roots are   
$$\Delta_{1}^+= \{\epsilon_i - \delta_1 | 1 \leq i \leq m \}.$$
There is only one simple root which is odd: $\epsilon_m-\delta_1$.\\

An element in $\frak h^*$ is called a weight. A weight $\Lambda$ will be denoted as follows:
\begin{equation}\label{eq:lamda1}
\Lambda 
= \sum_{i=1}^m \lambda_i\epsilon_i +  \mu \delta_1
=:  (\lambda_1, \cdots, \lambda_m; \mu).
\end{equation}
%
A weight $\Lambda$ is called integral if and only if $\lambda_i, \mu \in \mathbb{Z}$; it is called integral dominant if and only if it is integral and such that $\lambda_1 \geq \lambda_2 \geq \cdots \geq \lambda_m$.
As usual, we consider the following important weights.
$$\rho_0 =   \frac{1}{2} \sum_{\alpha \in \Delta _0^+} \alpha = \frac{1} {2}(m-1, m - 3, \ldots, 1 - m; 0),$$
$$  \rho_1 = \frac{1}{2}\sum_{\alpha \in \Delta_{1}^+}\alpha=\frac{1} {2}(1, 1, \ldots, 1 ; -m), \quad \rho = (m, \ldots, 2, 1; -1).
$$
There is a symmetric bilinear form $(\;,\;)$ on $\frak h^*$  is defined by 
$$(\epsilon_i, \epsilon_j) = \delta_{ij}, (\epsilon_i, \delta_1) = 0, (\delta_1, \delta_1) = -1.$$
The Weyl group of $\frak g $ is the Weyl group $ W$ of $\frak g_0$, hence it is the symmetric groups $ S_m $. For $ w \in W$, we denote by $\epsilon(w)$ its signature.\\
\subsection{Typical and atypical weights} 
Let $\Lambda = (\lambda_1, \cdots, \lambda_m; \mu)$ be a dominant weight.
A positive odd root $\epsilon_i - \delta_1 $, with $i = 1, 2, \ldots, m$, is said to be an atypical root of $\Lambda$ if
\begin{align}
(\Lambda + \rho, \epsilon_i - \delta_1) = 0.
\end{align}
Explicitly, this condition reads: $  \lambda_i +m + 1 - i =  - \mu + 1$ \cite{Kac2,Kac3}.
Denote by $\Gamma_\Lambda$ the set of atypical roots of $\Lambda$:
$$\Gamma_\Lambda = \{\epsilon_i - \delta_1| (\Lambda + \rho, \epsilon_i - \delta_1) = 0 \}.$$
Thus the number of elements of $\Gamma_\Lambda $ is $0$ or $1$. A weight $\Lambda$ is called typical if $\#\Gamma_\Lambda =0$ and atypical if  $\#\Gamma_\Lambda  = 1$.

\subsection{Kac modules}
For every integral dominant weight $\Lambda$, we denote by $V^0(\Lambda)$ the finite dimension irreducible $ \frak{g}_{\bar{0}}$-module with highest weight $\Lambda$, $V^0(\Lambda)$ is $(\frak{g}_{\bar{0}} \oplus \frak{g}_{+1})$- module with $\frak{g}_{+1}$ acting by $0$, where $\frak{g}_{+1}$ is the set of matries of the form $\left (\begin{array}{cc}0&B\\
0&0 \end{array}\right)$. Set
$$ 
\bar{V}(\Lambda) := \text{\rm Ind}_{\frak{g}_{\bar{0}} \oplus \frak{g}_{+1}}^{\frak{g}}V^0(\Lambda),
 $$
$\bar{V}(\Lambda)$ contains a unique maximal submodule $M(\Lambda)$. So, $\bar{V}(\Lambda)/M(\Lambda)$ is the irreducible module. Put
$$ V(\Lambda) := \bar{V}(\Lambda)/M(\Lambda).  $$
Then $V(\Lambda)$ is an irreducible module with highest weight $\Lambda$, it is called Verma module or Kac module \cite{Kac1}. 

\subsection{Characters of $\frak{gl}(m|1)$ and symmetric functions}
Let $V(\Lambda)$ be an irreducible representation with highest weight $\Lambda$ of $\frak g$. Such representations are $\frak h$-diagonalizable with weight decomposition $ V(\Lambda) = \bigoplus_{\mu}V_\mu$, where $V_\mu = \{ v \in V | h v = \mu(h) v \;\; \mbox{for all}\;\; h \in \frak{h}^*\}$ and character is defined to be 
$$\ch \; V = \sum_{\mu}(\dim V_\mu) e^{\mu},$$
 where $e^{\mu}$ ($\mu \in \frak h^*$) is the formal exponential. 

The character formula of irreducible representation of $\frak{gl}(m|1)$ below is the special case of character formula of irreducible representation of $\frak{gl(m|n)}$ is due to Su-Zhang \cite[Theorem 4.9]{Zhang2}.
\begin{equation}\label{eq:Su-Zhang}
 \ch V(\Lambda) = 
 \frac{1}{L_0}\sum_{w \in W}\epsilon (w)w \left( e^{\Lambda + \rho_0}\prod_{\beta \in \Delta_1^+ \backslash{\Gamma_\Lambda}}(1 + e^{-\beta}) \right), 
\end{equation}
where $L_0 = \prod_{\beta \in \Delta _0^+} (e^{\beta /2} - e^{-\beta /2}).$
We set $x_i:=e^{\epsilon_i}$ and $y:=e^{\delta_1}$. Then
\begin{equation}
L_0 = \prod_{\beta \in \Delta _0^+} (e^{\beta /2} - e^{-\beta /2}) =  \frac{\prod_{1 \leq i < j \leq m}(x_i - x_j)} {(\prod _{i = 1}^m x_i)^{(m - 1)/2}}.
\label{Eq:L0}
\end{equation}
and
 \begin{equation}\label{Eq:factor}
\prod_{\beta \in \Delta_1^+}(1+ e^{-\beta}) = \frac{\prod_{i = 1}^m(x_i + y)} {\prod_{i = 1}^m x_i }.
\end{equation}
Hence, if $\Gamma_\Lambda = \{\epsilon_k - \delta_1 \}$, we have
\begin{align*}
\prod_{\beta\in\Delta_1^+ \backslash{\Gamma_\Lambda} }(1 + e^{-\beta})  &= \frac{(x_1 + y)(x_2+y)\ldots (x_{k - 1} + y)(x_{k + 1}+y). . . (x_m + y))}{x_  1 x_2 \ldots x_{k - 1}x_{k + 1}  . . .x_m} \\
&= \frac{\sum_{i = 0}^{m - 1} e_i^{(k)}(x)}{x_1 x_2\ldots x_{k - 1}x_{k+1}\ldots x_m}y^{m -1 -  i},
\end{align*}
where $e_i^{(k)}(x)$ is the elementary symmetric function of $m - 1 $ variables $x_1, x_2, \ldots , x_{k - 1}$, $x_{k + 1},\ldots ,x_m.$
Recall that  
\begin{align*}
\rho_0 &= \frac{1}{2} \sum_{\alpha \in \Delta _0^+} \alpha = \frac{1} {2}(m-1, m - 3, \ldots, 1 - m).
\end{align*}
So, we have
\begin{equation}\label{Eq:exp_Lambda}
e^{\Lambda + \rho_0} = 
y^{\mu}    \prod_{i = 1}^m x_i^{\lambda_i + \frac{1}{2}(m- (2i - 1))}.
\end{equation}

For a partition  $\lambda = (\lambda_1,\ldots,\lambda_m)$ the corresponding {\em Schur symmetric function} in the variables $x_1,\ldots,x_m$ is defined as follows:
\begin{equation}\label{eq:38}
s_\lambda (x)= \frac{a_{\lambda +\delta}(x)}{a_{\delta}(x)},
\end{equation} 
where $\delta = (m-1, \ldots, 1, 0)$ and
\begin{equation}\label{eq:39}
 a _{\alpha }(x) := 
 \left |\begin{array}{cccc}x_{1}^{\alpha_1 }&x_1^{\alpha_2 }&\cdots&x_1^{\alpha_m }\\
x_{2}^{\alpha_1 }&x_2^{\alpha_2 }&\cdots&x_2^{\alpha_m}\\
\vdots&\vdots&\ddots&\vdots\\
x_{m}^{\alpha_1}&x_m^{\alpha_2}&\cdots&x_m^{\alpha_m }
\end{array}\right|  = \sum_{w \in S_m}\epsilon(w)
 w(x_1^{\alpha_1} x_2^{\alpha_2}\ldots x_m^{\alpha_m}),\end{equation}
where $w$ acts on the monomial 
$x_1^{\alpha_1} x_2^{\alpha_2}\ldots x_m^{\alpha_m}$ by permuting the variables (while keeping the exponents).
We shall extend the definition of $s_\lambda$ to any dominant weight by setting
$$s_\lambda:=\frac{a_{\lambda+\delta}}{a_\delta}.$$
With this notations $s_\lambda$ is exactly  the character of the irreducible 
representation of $\mathfrak{gl}(m)$ with highest weight equal to $\lambda$.
Note that for any integers $t$
$$s_\lambda= e_m(x)^ts_{(\lambda_1-t,\ldots,\lambda_m-t)}(x),$$
where $e_m(x):=\prod_ix_i$ -- the $m$-th elementary symmetric function.

\begin{lem}\label{L:2}
Let $\Lambda = (\lambda_1,\lambda_2,\ldots ,\lambda_m; \mu)$ be a weight. 

a) If $\Lambda$ is typical then
\begin{equation}
 \ch V(\Lambda) = \frac{y^{\mu}}{ e_m(x) }  s_{\lambda }(x) \sum_{i = 0}^m e_i(x) y^{m - i},
 \end{equation}
where $\lambda= (\lambda_1,\lambda_2,\ldots ,\lambda_m)$, $e_1(x),e_2(x), \ldots , e_m(x)$ are elementary symmetric 
functions of in ${x_1, x_2, \ldots, x_m}$.\\
b) If $\Lambda $ is 
atypical with $\Gamma_\Lambda = \{\epsilon_k - \delta_1\}$, for some $k \in \{ 1, 2, \ldots, m\}$, then 
\begin{equation} \ch V(\Lambda)  = \frac{ y^{ \mu }}{ e_m(x)}  
\sum_{i = 1}^{ m} \left(\sum_{\alpha\in A_{k,i}} s_{\alpha}(x) \right)y^{m -(i+1)},
\end{equation}
where $A_{k,i}$ is the set of integral dominant weights $\alpha$ 
(of size $m$) such that $\alpha_j-\lambda_j\in\{0,1\}$ and $|\alpha|-|\lambda|=i + 1$, furthermore, $\alpha_k-\lambda_k=1$.
\end{lem}
\begin{proof}[Proof]
The case $\Lambda$ being typical can be readily seen from Equations \eqref{Eq:L0}, \eqref{Eq:factor}, \eqref{Eq:exp_Lambda}.
Assume now  $\Lambda = (\lambda_1,\lambda_2,\ldots ,\lambda_m; \mu)$ atypical with $\Gamma_\Lambda = \{\epsilon_k - \delta_1\}$ for some $k$, $1\leq k\leq m$. Set
$$ D:= \prod_{1 \leq i< j \leq m}(x_i - x_j).$$
Following \eqref{eq:Su-Zhang} we have 
\begin{align*}
\ch V(\Lambda)&= \frac{1}{L_0} \sum_{w \in S_m} \epsilon(w) w \left(e^{(\Lambda + \rho_0 )}\prod_{\beta\in\Delta_1^+ \backslash{\Gamma_\Lambda} }(1 + e^{-\beta})\right)\\ 
  &= \frac{y^{\mu}}{ e_m(x) D  } 
    \sum_{w \in S_m} \epsilon(w)w \left( \prod_{i = 1}^m x_i^{\lambda_i + (m- i)} x_k  \sum_{i = 0}^{ m - 1}e_i^{(k)}(x) y^{m -1 -  i}\right)\\ 
  & = \frac{y^{\mu} }{e_m (x)  D  } \sum_{i = 0}^{ m - 1} \left(\sum_{j\in J_i} a_{\lambda+j+\delta} (x)  \right)y^{m -1 -  i}\\
&=\frac{y^{ \mu }}{ e_m (x)} \sum_{i = 0}^{ m - 1} \left(\sum_{\alpha\in A_{k,i}} s_{\alpha}  (x)\right)y^{m -1 -  i},
\end{align*}
where $J_i$ denotes the set of $m$-tuples $(j_1 ,j_2 ,\ldots,j_m )$, in which $i+1$ entries are equal to  1 and the others are equal to $0$, further $ j_{k}  = 1$.
\end{proof}
\section{Composite partitions and Schur functions}
 Our aim is to give a character formula for irreducible representations of $\mathfrak g$, the weights of which belong to a class of integral dominant weights specified below. We notice that an arbitrary integral dominant weight is different from a weight in this class by a multiple of the weight $\sigma:=(1,...,1;-1)$ which corresponds to the super-determinantal representation (see Lemma \ref{L:3}).
 
  An integral dominant weight $$ \Lambda = (\lambda_1, \lambda_2, \ldots, \lambda_m; -k), $$ 
with $ 0 \leq k \leq m$ and $\lambda_{m-k}\geq 0 \geq \lambda_{m-k+1}$, is called a {\em special weight}.
 We denote by $P$ the set of all special weights and for each $k$, $0 \leq k \leq m$, let $P_k$ be the subset:
\begin{equation}\label{eq:15}
P_k = \{ \Lambda = (\lambda_1, \lambda_2, \ldots, \lambda_m; -k),  | \lambda_{m-k}\geq 0 \geq \lambda_{m-k+1}     \}.
\end{equation} 
 In particular, $P_0$ is the set of genuine partitions, i.e. weights of covariant irreducible representations  (i.e. those irreducible which can be constructed by multi-linear algebra from the fundamental representation).

 Our first task is to transform a special weight into an $m$-standard composite weight. Then we define a super-symmetric function in term of the $m$-standard partition using a Jacobi-Trudi type formula. This super-symmetric function turns out to be the irreducible character corresponding to the original weight.

\subsection{Composite partitions}
A composite partition is a pair of partition $(\bar\nu;\mu)$, it is called {\em $m$-standard} if
\begin{equation}\label{eq:42}
l(\mu) + l(\nu) \leq m.
\end{equation}
Let $Q$ be the set of $m$-standard composite partitions and for each $ 0 \leq k \leq m$
let $Q_k$ be the subset of those $(\nu;\mu)$ with $l(\nu) =k$:
\begin{equation}\label{eq:}
 Q_k = \{ (\nu;\mu) | l(\mu) \leq m - k , l(\nu) = k\}.
\end{equation} 

We define a map $\varphi:P\to Q$ as follows. For $\Lambda\in P$, $\Lambda = (\lambda_1, \lambda_2, \ldots, \lambda_m; -k)$, set
 $\varphi(\Lambda)$ to be the composite partition $(\nu,\mu)$, where
 $$
\begin{cases}
    \mu &= (\lambda_1, \lambda_2, \ldots, \lambda_{m - k}),   \\ 
    \nu &= (1 - \lambda_m, 1 - \lambda_{m - 1}, \ldots, 1 - \lambda_{m - k+1}).
\end{cases} 
$$
Notice that $l(\mu)\leq m-k$ and $l(\nu)=k$. Thus, if $\Lambda\in P_k$, 
$( \nu;\mu) \in Q_k.$  

\begin{lem}\label{lemPQ}
The map $\varphi$ is bijective and restricts to a bijective map from $P_k$ to $Q_k$ for each $0\leq k\leq m$. Further, 
$\Lambda$ is typical iff $l(\mu)+l(\nu)=m$. 
\end{lem}
\begin{proof}[Proof]
It is easy to see that $\varphi$ is injective and maps $P_k$ to $Q_k$.  
To describe the converse map we shall deploy the following notation, which will be used later.
Let $(\nu,\mu)$ be an $m$-standard composite partition. 
We define its $m$-composition as a weight of the following form:
\begin{equation}\label{Eq:cup_op}
(\mu\cup_m\bar\nu):=(\mu_1,\ldots,\mu_{m-k},-\nu_k,\ldots,-\nu_1).\end{equation}
Now  define $\Lambda$ to be
$$(\mu \cup_m((1^{l(\nu)})+\bar\nu);-l(\nu)).$$ 
Then $\Lambda$ is the pre-image of $(\nu,\mu)$ under $\varphi$. Explicitly,
$$\Lambda=(\mu_1,\ldots,\mu_{m-k},1-\nu_k,\ldots,1-\nu_1;-l(\nu) ).$$

The last claim, in the above notations, amounts to saying that $\lambda$ is typical iff $\mu_{m-k}\neq 0$.
If $\mu_{m - k} = 0$ then $\Gamma_\Lambda = \{\epsilon_{m - k} - \delta_1 \}$. So $\Lambda$ is atypical.

Conversely, if $\mu_{m - k} \ne 0$ we show that $\Gamma_\Lambda =  \emptyset$, which amounts to proving
 $  \lambda_i +m + 1 - i \ne  k + 1, \forall i \in \{ 1,2,\ldots, m\} $. 
 
 Indeed, for $i \leq m - k$ we have $\lambda_i = \mu_i > 0$ and
$$ \lambda_i + m + 1 - i = \mu_i + m + 1 -i \geq \mu_i + m + 1 -( m -k) > k + 1.
 $$
For $i > m - k$ then $\lambda_i \leq 0$ and
$$ \lambda_i + m + 1 - i < \lambda_i + m + 1 - ( m -k) \leq  k + 1.
 $$
Thus $\Lambda$ is typical.
\end{proof}  

\subsection{Symmetric functions associated to composite partitions}
Let $(\nu; \mu)$ be an $m$-standard composite partition  (cf. \eqref{eq:42}). Set   $l(\mu)= p, l(\nu) = q$.  The {\em symmetric function indexed by this composite partition} is defined by:
\begin{equation}\label{eq:51}
s_{(\nu; \mu)}(x) :=s_{(\mu\cup_m\bar\nu)}(x),
\end{equation} 
where as defined in the previous paragraph.
$$(\mu\cup_m\bar\nu)= (\mu_1, \mu_2, \ldots, \mu_p,0,\ldots,0,-\nu_q ,\ldots,-\nu_2,-  \nu_1).$$
The following formula for symmetric functions indexed by 
composite partitions in terms of elementary symmetric functions
 was conjectured by Balantekin and Bars \cite{bar2}, and proved in \cite{CK}, namely
\begin{equation}\label{eq:52}
s_{(\nu;\mu)} (x)= \det \left(
\begin{tabular}{c|c}
$\dot{h}_{\nu_l+k-l}(x)$&$h_{\mu_j-k-j+1}(x)$\\
\hline
$\dot{h}_{\nu_l-i-l+1}(x)$&$h_{\mu_j +i -j}(x)$
\end{tabular}\right),
\end{equation}
where the indices $i, j, k$ resp. $l$ run from top to bottom, from left to right, from bottom to top resp. from right to left and with the function $\dot{h_r}(x) = h_r(x_1^{-1}, \ldots, x_m^{-1})$.\\

\subsection{Super-symmetric functions associated to composite partitions}
Consider two sets of independent variables $x = \{x_1, x_2, \ldots, x_m\}$ and $\{y\}$. The {\em complete super-symmetric} functions can be written in terms of the elementary symmetric and the complete symmetric functions:
\begin{equation}\label{eq:}
h_r(x/y) =\sum_{k = 0}^r h_k(x) e_{r-k}(y)= \sum_{k = 0}^r h_k(x) y^{r-k}.
\end{equation} 
Given any partition $\lambda = (\lambda_1, \lambda_2, \ldots)$, one defines  the corresponding {\em super-symmetric Schur fuctions} to be 
\begin{equation}\label{eq:53}
s_\lambda(x/y) = \det (h_{\lambda_i - i + j}(x/y))_{1\leq i,j \leq l(\lambda)}.
\end{equation} 
In particular, $s_{(r)}(x/y) = h_r(x/y).$ 
Given a composite partition $(\nu; \mu)$, one defines  the associated  {\em super-symmetric function} to be (\cite{bar1,bar2})
\begin{equation}\label{eq:54}
 s_{({\nu};\mu)} (x/y)= \det \left(
\begin{tabular}{c|c}
$\dot{h}_{\nu_l+k-l}(x/y)$&$h_{\mu_j-k-j+1}(x/y)$\\
\hline
$\dot{h}_{\nu_l-i-l+1}(x/y)$&$h_{\mu_j +i -j}(x/y)$
\end{tabular}\right),
\end{equation}
where the indices $i$ resp. $j,k, l$ run from top to bottom, resp. from left to right, from bottom to top, from right to left and $\dot{h_r}(x/y) = h_r(\bar{x}/\bar{y})$ with $\bar{x_i} = x_i^{-1}$, $\bar{y} = y^{-1}$. 


\begin{lem}\label{L:51}
 Let $({\nu}; \mu)$ be a composite partition. Then
$$ s_{({\nu}; \mu)}(x/y) = \sum_{\alpha, \beta}s_{(\beta; \alpha)}(x) y^{a-b},$$
where $a= \left|\mu -\alpha \right|, b= \left|\nu -\beta \right|$ and the sum is taken over all partitions $\alpha$ and $\beta$ such that $ (\mu - \alpha)_i, (\nu - \beta)_i \in \{0,1\}$.
\end{lem}
\begin{proof}
This follows from \cite[Lemma A.3]{M4}.
\end{proof}

\section {The main theorem}
 

\begin{thm}\label{T:1}
Let $ \Lambda = (\lambda_1, \lambda_2, \ldots, \lambda_m; -k)$ be a special weight in $P$ and let  $(\nu;\mu) \in Q$ be the corresponding composite partition. Then
$$  {\rm ch} V(\Lambda) = s_{(\nu;\mu)}(x/y). $$
\end{thm} 
\begin{proof}
 We consider two cases: $\Lambda$ is typical and $\Lambda$ is atypical. Recall that we have the following relationship between $\Lambda$ and $\mu,\nu$:
 $$\Lambda=(\lambda;-k)=(\mu\cup_m((1^k)+\bar\nu))=(\mu_1,\ldots,\mu_{m-k},1-\nu_k,\ldots,1-\nu_1).$$

Case I:  $\Lambda$ is typical, thus $\mu_{m-k}\neq 0$.
According to Lemmas  \ref{L:2} and  \ref{L:51} we need to show that
\begin{equation}\label{Eq:typical}  s_{\lambda }(x) \sum_{i = 0}^m e_i(x) y^{m - i-k} =  \sum_{\alpha, \beta}   e_m(x)s_{(\beta;\alpha)}(x)y^{a-b},
 \end{equation}
where, on the right-hand side, $a= \left|\mu -\alpha \right|, b= \left|\nu -\beta \right|$ and $(\mu - \alpha)_i, (\nu - \beta)_i \in \{0,1\}.$\\
This yields,\\
$$\left\{ \begin{array}{l}
\alpha = (\mu_1- j_1^\alpha,\ldots, \mu_{m - k}-j_{m-k}^\alpha),
\\
\beta =(\nu_1- j_1^\beta, \ldots, \nu_k-j_k^\beta),
\end{array} \right.
$$
where $(j^\alpha_1, \ldots, j^\alpha_{m-k})$ only includes 0 and 1, it has $a$ numbers 1 and $m-k-a$ numbers 0; similarly, $(j^\beta_1, \ldots, j^\beta_k)$ has $b$ numbers 1 and $k-b$ numbers 0.\\
By \eqref{eq:51}, the right-hand side of \eqref{Eq:typical} is equal to
\begin{align*}
 e_m(x) \sum_{(j^\alpha_1, \ldots, j^\alpha_{m-k})}\sum_{(j^\beta_1, \ldots, j^\beta_k)}s_{(\mu_1-j_1^\alpha,\ldots,\mu_{m- k}-j_{m-k}^\alpha,   -\nu_k +j_k^\beta,\ldots, -\nu_1+j_1^\beta)}(x) y^{a-b}.
 \end{align*}
Notice that $e_m(x)s_{\lambda}(x)=s_{\lambda+(1^m)}(x)$. Hence the coefficient of $y^{m-i-k}, 0\leq i \leq m$, in this sum is:
\begin{align*}
 C_{m-i-k}= \sum_{(j^\alpha_1, \ldots, j^\alpha_{m-k})}\sum_{(j^\beta_1, \ldots, j^\beta_k)}
s_{(\mu_1+j_1^\alpha,\ldots,\mu_{m-k}+j_{m-k}^\alpha,1-\nu_k+j_k^\beta,\ldots, 1-\nu_1+j_1^\beta)}(x),
\end{align*}
where the sums run on tuples: $(j^\alpha_1, \ldots, j^\alpha_{m-k}) \; \mbox{and} \; (j^\beta_1, \ldots, j^\beta_k)$ such that the number of 1's in $(j^\alpha_1, \ldots, j^\alpha_{m-k})$ {  plus} the number of 1's in $(j^\beta_1, \ldots, j^\beta_k)$ is $i$.

\begin{align*}
C_{ m-i-k} = \sum_{J} 
 s_{(\mu\cup_m(1^k)+\bar\nu)+J}(x),
 \end{align*}
 where $J$ runs in the set of $m$-tuples $(j_1,\ldots,j_m)$ with $j_i\in\{0,1\}$ and $|J|=i$.
According to the Robinson-Schensted rule we conclude that
$$C_{m-i-k}=s_{\lambda}(x)e_i(x).$$
This finishes the proof of \eqref{Eq:typical}.

Case II: $\Lambda$ is atypical, thus $\mu_{m-k}=0$.
According to Lemmas  \ref{L:2} and  \ref{L:51} we need to show that
\begin{equation}\label{Eq:atypical}    
\sum_{i = 1}^{ m} \left(\sum_{\alpha\in A_{m-k,i}} s_{\alpha}(x) \right)y^{m -(i+1)-k} =  \sum_{\alpha, \beta}   e_m(x)s_{(\beta;\alpha)}(x)y^{a-b},
 \end{equation}
where, on the right-hand side, $a= \left|\mu -\alpha \right|, b= \left|\nu -\beta \right|$ and $(\mu - \alpha)_i, (\nu - \beta)_i \in \{0,1\}.$\\
The last conditions amount to
$$\left\{ \begin{array}{l}
\alpha = (\mu_1- j_1^\alpha,\ldots, \mu_{m - (k+1)}-j_{m-(k+1)}^\alpha),
\\
\beta =(\nu_1- j_1^\beta, \ldots, \nu_k-j_k^\beta),
\end{array} \right.
$$
where the sequence $(j^\alpha_1, \ldots, j^\alpha_{m-(k+1)})$ contains $a$ times 1 and $m-(k+1)-a$ times  0; similary, $(j^\beta_1, \ldots, j^\beta_k)$ contains  $b$ times 1 and $k-b$ times zeros.\\
By \eqref{eq:51}, the right-hand side of \eqref{Eq:atypical} is equal to
$$ 
e_m(x) \sum_{(j^\alpha_1, \ldots, j^\alpha_{m-k-1})}\sum_{(j^\beta_1, \ldots, j^\beta_k)}s_{(\mu_1 -j_1^\alpha,\ldots,\mu_{m-(k+1)}-j_{m-k-1}^\alpha, 0,  -\nu_k+j_k^\beta,\ldots,-\nu_2+ j_2^\beta,-\nu_1+j_1^\beta)}(x) y^{a-b}.
 $$
Hence the coefficient of $y^{m-(i+1)-k}, 0\leq i \leq m - 1$, in this sum is:
\begin{eqnarray*}
 \lefteqn{C_{m-(i+1)-k}=}\hspace{\textwidth}  \\
  \sum_{(j^\alpha_1, \ldots, j^\alpha_{m-k-1})}\sum_{(j^\beta_1, \ldots, j^\beta_k)}s_{(\mu_1+(1-j_1^\alpha) , \mu_2+(1-j_2^\alpha),\ldots,\mu_{m-(k+1)}+ (1 - j_{m-(k+1)}^\alpha), 1, 1 -\nu_k+ j_k^\beta,\ldots,1 -\nu_1+ j_1^\beta)}(x),
\end{eqnarray*} 
where the sums run on tuples: $(j^\alpha_1, \ldots, j^\alpha_{m-(k+1)}), (j^\beta_1, \ldots, j^\beta_k)$ such that the number of 1's in $(j^\alpha_1, \ldots, j^\alpha_{m-(k+1)})$ minus the number of 1's in $(j^\beta_1, \ldots, j^\beta_k)$ equals $m-(i + 1)-k$.\\
Set\\
$$\left\{ \begin{array}{l}
r_l = 1 - j_l^{\alpha} \quad 1\leq l \leq m-k -1\;\; \mbox{and} \;\; r_{m-k} = 1
\\
r_{m-k + l} = j_{k -l+1}^{\beta} \quad 1 \leq l \leq k.
\end{array} \right.$$
Then
\begin{eqnarray*}
\lefteqn{ C_{m-(i+1)-k}= }\hspace{\textwidth}\\
\sum_{(r_1, r_2, \ldots, r_{m})}s_{(\mu_1+ r_1 , \mu_2+ r_2,\ldots,\mu_{m-(k+1)}+r_{m-(k+1))}, r_{m-k}, 1 -\nu_k+ r_{m-k+1},\ldots,1 -\nu_2+ r_2,1 -\nu_1+r_1)}(x),
\end{eqnarray*}
where $(r_1, r_2, \ldots, r_{m})$ is $m$-tuple, in which $i+1$ entries are equal to  1 and the others are equal to $0$, further $ r_{m-k}  = 1$. Consequently 
$$
 C_{m-(i+1)-k}= \sum_{\eta\in A_{m-k,i}} s_{\eta}(x),
 $$
where $A_{m-k,i}$ is the set of integral dominant weights $\eta$ 
(of size $m$) such that $\eta_j-\lambda_j\in\{0,1\}$ and $|\eta|-|\lambda|=i + 1$, furthermore, $\eta_{m-k}-\lambda_{m-k}=1$.
This finishes the proof of \eqref{Eq:atypical}.
\end{proof}

\subsection{Special  weights} We have so far established a determinantal formula of Jacobi-Trudi type for the class of special weights. As in the classical case, we show that this essentially furnishes all irreducible characters of $\mathfrak{gl}(m|1)$. Namely, we show that an integral dominant weight $\Lambda$ can be represented in the form 
$$\Lambda=\Lambda'+j\sigma,$$
where $\Lambda'$ is a special weight and 
$$\sigma=(1,1,\ldots,1;-1).$$
The following formula is well-known:
\begin{equation}\label{eq:char}
\ch \;V(\Lambda+j\sigma) = (e^\sigma)^j \ch \;V(\Lambda),
\end{equation}  
where $e^\sigma$ the formal exponential, $e^\sigma=y^{-1}\prod_ix_i$.
\begin{pro}\label{L:3}
Let $\Lambda $ be an integral dominant weight. Then there is unique integer $j$ such that $\Lambda' :=\Lambda + j\sigma $ has the following form:
 $$ \Lambda' = (\lambda_1, \lambda_2, \ldots, \lambda_m; -k), $$ 
with $ 0 \leq k \leq m$ and $\lambda_{m-k}\geq 0 \geq \lambda_{m-k+1}$, where $\sigma = (1, \ldots, 1; -1)$. That is $\Lambda'$ is a special weight.
\end{pro}
\begin{proof}[Proof]
We use induction on $m$. The case $m=1$ is easy, a weight $(\lambda;-\mu)$ is led to 
$$\left[\begin{array}{ll} (\lambda - \mu;0)& \text{ if } \lambda\geq \mu\\
(\lambda-\mu+1;-1)&\text{ if } \lambda<\mu.\end{array}\right.$$

Let's consider the induction step $m$ to $m+1$. Given a weight
$$(\alpha_1,\alpha_2,\ldots,\alpha_{m+1};-\beta),$$
 using the induction hypothesis on the weight
$$(\alpha_2,\ldots,\alpha_{m+1};-\beta),$$
we can bring it to the form 
 $$(\lambda_1,\lambda_2,\ldots,\lambda_{m+1};-k),$$
with $0\leq k\leq m$ (by adding a multiple of $\sigma$), such that
$$\lambda_{(m+1)-k}\geq 0\geq\lambda_{(m+1)-k+1}.$$
Note the shift of the indexes and the condition $\lambda_1\geq 0$ is not imposed (when $k=m$).

Thus, if in this new weight we have $k<m$ then it automatically satisfies the requirement. Similarly, if in this new weight we have $k=m$ and $\lambda_1\geq 0$ then it also satisfies.

 It remains the case $k=m$ and $\lambda_1<0$. Then adding (or substracting) a multiple of $\sigma$ to this weight we get the partition
$$(\lambda_1+1,\lambda_2+1,\ldots,\lambda_{m+1}+1;-(m+1))$$
satisfying $0\geq\lambda_1 + 1$.
 
Finally, we prove the uniqueness assertion. Let $\Lambda$ be a integral dominant weight. Assuming $j, j'$ are integers such that\\
$$
\begin{cases}
\Lambda + j \sigma = (\lambda_1, \ldots, \lambda_m; -k) \quad \mbox{with}\;  \lambda_{m-k}\geq 0 \geq \lambda_{m-k+1}\\
\Lambda + j' \sigma = (\lambda'_1, \ldots, \lambda'_m; -k') \quad \mbox{with}\;  \lambda'_{m-k'}\geq 0 \geq \lambda'_{m-k' + 1} .
\end{cases}
$$
We need show that $j = j'$. Assume the contrary, then we can assume  $ j > j'$. We have\\
$$ (j - j')\sigma = \Lambda + j \sigma - (\Lambda + j'\sigma) =  (\lambda_1 - \lambda'_1, \ldots, \lambda_m - \lambda'_m; - (k -  k')).  $$
Then
$$  j - j' = k - k' := t .$$
This yields
$$ t > 0,\; k = k' + t \;\mbox{and} \; j = j' + t .$$
On the other hand, 
$$ \Lambda + j\sigma = \Lambda + j' \sigma + t\sigma = (\lambda'_1 + t, \ldots, \lambda'_m + t; -(k' + t)). $$
Since
$$   \lambda_{m-k}\geq 0 \geq \lambda_{m-k+1},$$
we have
$$  \lambda'_{m- k}+ t \geq 0 \geq \lambda'_{m- k +1} + t .$$
This implies $$0 \geq \lambda'_{m-(k' + t)+1} + t \geq \lambda'_{m- k'} + t \geq t > 0,$$
which is a contradiction. Thus we conclude $j = j'$.
\end{proof}
 \begin{cor}
 The characteristic of irreducible representations of $\frak{gl}(m|1)$ are represented by super-symmetric S-functions.
\end{cor}

\section{Concluding remarks} 
The starting point of this work is the construction of irreducible representations of $\mathfrak{gl}(3|1)$, given in \cite{dhh,dung,dh}, which gives a determinantal type formula expressing an irreducible representations in terms of the symmetric powers of the fundamental representation and their duals. We later discovered that a vast generalization of this formula has been provided by Moens and van der Jeugt.
In the paper \cite[2004]{M4}, E.M. Moens and J. van der Jeugt announce the following theorem.\\

\noindent{\bf Theorem.} \cite[Theorem 4.3]{M4}  {\em
Let $(\bar{\nu}; \mu)$ be a standard and critical composite partition (see \cite[Definition 3.1]{M4}) with no overlap (see \cite[Section 2]{M4}) and $\Lambda_{(\bar{\nu}; \mu)}$ be the corresponding super weight. The character $  {\rm ch} V(\Lambda_{(\bar{\nu}; \mu)}) $ is equal to $s_{(\bar{\nu};\mu)}(x/y) $, which is defined in the same manner as in \eqref{eq:54}.
 }\\

The proof of this theorem is based on the following lemma. 

\noindent{\bf Lemma.}  \cite[Lemma A.5]{M4} {\em 
Suppose $|x| = m, |y| = n$ and $h, p$, $q$ are positive integers with $ m = p + q$. Let $\kappa=(\kappa_1, \kappa_2, \ldots, \kappa_q), \eta =(\eta_1, \eta_2, \ldots)$ and $\mu$ be partitions, and $\nu = (\kappa_1, \kappa_2, \ldots, \kappa_q, \eta_1, \eta_2, \ldots)$. Then
\begin{equation*}
\sum_{x' + x''}\frac{(\prod x')^q (\prod x'')^h s_{(\overline\eta; \mu)}(x'/y) s_{\kappa + (h^q)}( \bar{x''}/ \bar{y})}{E(x', x'')} = s_{(\overline \nu; \mu)}(x/y),
\end{equation*}
where the sum is over all possible decompositions $x = x' + x''$ with $|x'| = p $ , $|x''| = q$.
}\\

However, Moens notices in his thesis   that this Lemma is false and proposes to prove the above theorem by using a weaker form of this lemma Lemma, in which one assumes that $\bar{\nu}; \mu$ is a critical composite partition with no zeros in the overlap when presented in the $m\times n$-rectangle \cite[Lemma 5.14]{M6}. 
However no proofs are provided.
 Thus the mentioned above theorem  in its general form is still  a conjecture. 

\section{Acknowledgment}
This research  is funded by  Vietnam National Foundation for Science and Technology Development (NAFOSTED), grant number 101.04-2016.19. A part of this work was carried out when the first and the third authors were visiting the Vietnam Institute for Advanced Study in Mathematics. They would like to thank VIASM for the financial support and the excellent working environment.

\end{document}